\documentclass[12pt]{article}  

\usepackage{amsmath,amsfonts,latexsym,amsthm,amssymb,mathabx}

\title{Fractional approximation of solutions of evolution equations}
\author{Anatoly N. Kochubei, Yuri G. Kondratiev}
\date{}
\newcommand{\D}{\mathbb D_t^{(\alpha )}}

\newtheorem*{teo}{Theorem}
\newtheorem*{defin}{Definition}
\begin{document}

\maketitle\thispagestyle{empty}


\begin{abstract}
We show how to approximate a solution of the first order linear evolution equation, together with its possible analytic continuation, using a solution of the time-fractional equation of order $\delta >1$, where $\delta \to 1+0$.
\end{abstract}

\renewcommand{\thefootnote}{}
\footnotetext{\hspace*{-.51cm}AMS 2010 subject classification: Primary: 34G10,
35R11; Secondary: 30B40\\ %
}

\section{Introduction}\label{section1}

In many applications (see, for example, \cite{BKK,FKO}) we encounter evolution equations, which are too singular for application of the standard semigroup approach. A possible remedy is to consider the equation on a scale of Banach spaces and to use an approach known as Ovsyannikov's method; see \cite{F} for its history and complete references.

Barkova and Zabreiko \cite{BZ} extended this method to evolution equations with the Caputo-Djrbashian fractional derivative $\D$ of order $\alpha >0$. In particular, they showed that the Cauchy problem considered in this setting becomes less singular with the growth of $\alpha$. For example, it can happen that a solution $u$ of the first order equation $u_t'=Au$ exists only on a finite time interval, while the existence of solutions of the equation $\mathbb D_t^{\delta}u=Au$, $1<\delta <2$ can be guaranteed for all $t>0$. Therefore such solutions $u_\delta$ are natural means of approximating solutions of the initial first order equation.

In this note we prove an even stronger result -- the ``fractional approximations'' $u_\delta (t^{1/\delta })$ approach, as $\delta \to 1+0$, not only the solution $u$, but its maximal possible analytic continuation. As we will show, this follows from Hardy's theorem \cite{Ha} about approximate analytic continuation obtained by summing certain divergent series.

\bigskip

\section{Cauchy problems}\label{section2}

Consider the Cauchy problem
$$
\frac{du}{dt}=Au,\quad u(0)=u_0,
$$
in a scale of Banach spaces $X_\omega$, $\omega \in [0,1)$, $X_\omega'\subset X_\omega''$ for $\omega'<\omega''$. Here $A$ is a linear operator, which is bounded from $X_\omega'$ to $X_\omega''$ for each couple of indices with $\omega'<\omega''$, and
$$
\| A\|_{\omega'\to \omega''}\le \frac{C}{\omega''-\omega'}.
$$
For each couple $(\omega',\omega'')$ and any initial vector $u_0\in X_{\omega'}$, such that $Au_0\in X_{\omega'}$, there exists a local solution in $X_{\omega''}$ of the form
$$
u(t)=\sum \limits_{n=0}^\infty \frac{t^n}{n!}A^n u_0,\quad 0\le t<T_{\omega',\omega''};
$$
see \cite{BKK,Z}.

Let $1<\delta <2$. Consider also the Cauchy problem with the Caputo-Djrbashian fractional derivative
$$
\mathbb D_t^{(\delta )}u_\delta =Au_\delta ,\quad u_\delta (0)=u_0, u_\delta'(0)=0.
$$
Under the same assumptions, it has the solution (see \cite{BZ})
$$
u_\delta (t)=\sum \limits_{n=0}^\infty \frac{t^{\delta n}}{\Gamma (\delta n+1)}A^n u_0,
$$
which exists for all $t>0$ (the above series equivalent to the iteration process from \cite{BZ} converges on any finite interval).

In fact, the solution $u(t)$ can be continued to a holomorphic function on the disk $\{ t\in \mathbb C:\ |t|<T_{\omega',\omega''}\}$ while the function
$$
u_\delta (t^{1/\delta })=\sum \limits_{n=0}^\infty \frac{t^n}{\Gamma (\delta n+1)}A^n u_0
$$
is extended to an entire function. Both are with values in $X_{\omega''}$. Note that
$$
u_\delta (t^{1/\delta })=\sum \limits_{n=0}^\infty \lambda_n(\delta )\left( A^n u_0\right) \frac{t^n}{n!}
$$
where $\lambda_n(\delta )=\dfrac{n!}{\Gamma (\delta n+1)}$.

\bigskip

\section{Analytic continuations}\label{section3}

\begin{defin}[see \cite{Ha}]
Let $f(z)$ be a holomorphic function on a neighborhood of the point $z=0$ determined there by a convergent power series. The {\it Mittag-Leffler star} $G(f)$ of the function $f$ is a domain obtained from $\mathbb C$ as follows: draw a ray from the origin to each singular point of the function $f$, and cut the plane along the part of the ray located after the singular point.
\end{defin}

\medskip
For example, the star of the function $\sum\limits_{n=0}^\infty z^n=\dfrac1{1-z}$ is $\mathbb C\setminus [1,\infty)$.

\medskip
\begin{teo}
Suppose that the solution $u(t)$ is continued along rays to a single-valued holomorphic function on $G(u)$. Then
$$
u_\delta (t^{1/\delta })\longrightarrow u(t),\quad t\in G(u),
$$
uniformly on any closed bounded domain inside $G(u)$.
\end{teo}

\medskip
\begin{proof} By Theorem 135 (Section 8.10) from \cite{Ha} (evidently valid also for vector-functions), it is sufficient to prove that the function
$$
\varphi_\delta (z)=\sum \limits_{n=0}^\infty \dfrac{n!}{\Gamma (\delta n+1)}z^n,\quad \delta >1,
$$
is entire and tends to $\dfrac1{1-z}$ uniformly, as $\delta\to 1+0$, in each closed bounded domain not intersecting the semi-axis $[1,\infty )$. The entireness property follows from the Stirling formula. The above convergence is obvious for $|z|<1$. Due to the uniqueness theorem for holomorphic functions, it will be sufficient to prove that $\varphi_\delta (z)$ tends to some holomorphic function uniformly on compact subsets of $\mathbb C\setminus [0,\infty )$.

The function $\varphi_\delta$ can be represented as a Wright function:
$$
\varphi_\delta (z)={}_2\Psi_1\Bigl[ \begin{matrix}
(1,1) & (1,1) \\ (1,\delta )\end{matrix}\Bigl| z \Bigr] =\sum \limits_{k=0}^\infty \frac{[\Gamma (1+k)]^2}{\Gamma (1+\delta k)}\frac{z^k}{k!}
$$
with parameter $\Delta =\delta -2$, which means that $-1<\Delta <1$ (for information about the Wright functions see \cite{KST}, 1.11).

Let us use the integral representation for the Wright functions. In our case we have
$$
\varphi_\delta (z)=\frac1{2\pi i}\int\limits_{\frac12-i\infty}^{\frac12+i\infty}\frac{\Gamma (s)[\Gamma (1-s)]^2}{\Gamma (1-\delta s)}(-z)^{-s}ds.
$$
Since $\Gamma (s)\Gamma (1-s)=\frac{\pi}{\sin (\pi s)}$, we find that
\begin{equation}
\varphi_\delta (z)=\frac1{2i}\int\limits_{\frac12-i\infty}^{\frac12+i\infty}\frac{\Gamma (1-s)}{\sin (\pi s)\Gamma (1-\delta s)}(-z)^{-s}ds,\quad |\arg (-z)|<\frac{(3-\delta )\pi}2
\end{equation}
(see \cite{KST}, (1.11.21), (1.11.22)). Note that for $\delta$ close to 1 this representation becomes valid for any $z$ outside the semi-axis $(0,\infty )$.

By a corollary of the Stirling formula (\cite{KST}, (1.5.14)), for $s=\frac12+iy$, $|y|\to \infty$, we have
$$
\left| \frac{\Gamma (1-s)}{\Gamma (1-\delta s)}\right| \le C|y|^{\frac{\delta -1}2}\exp \left\{ \frac{\pi (\delta -1)|y|}2\right\}.
$$
Noticing also that $|\sin (\pi (\frac12 +iy)|=|\cos (i\pi y)|=\frac12 \left| e^{\pi y}+e^{-\pi y}\right|$ we see the possibility to pass to the limit in (3.1), as $\delta \to 1+0$. Thus, for $z\notin [0,\infty )$, uniformly on compact sets,
$$
\varphi_\delta (z) \longrightarrow  \frac1{2i}\int\limits_{\frac12-i\infty}^{\frac12+i\infty}\frac{(-z)^{-s}}{\sin (\pi s)}\,ds.
$$
The right-hand side is a holomorphic function on $\mathbb C\setminus [0,\infty )$. \end{proof}

\bigskip
\noindent
\parbox[t]{.48\textwidth}{
Anatoly N. Kochubei\\
Institute of Mathematics\\
National Academy of Sciences\\ of Ukraine \\
Tereshchenkivska 3\\
01601 Kyiv, Ukraine\\
kochubei@i.com.ua \\
kochubei@imath.kiev.ua} \hfill
\parbox[t]{.48\textwidth}{
Yuri G. Kondratiev\\
Fakult\"{a}t f\"{u}r Mathematik\\
Universit\"{a}t Bielefeld\\
33615 Bielefeld,Germany\\
kondrat@math.uni-bielefeld.de
}

\end{document}